\documentclass[11pt]{amsart}

\usepackage{amssymb}
\usepackage{amsthm}
\usepackage{amsmath}
\usepackage{amsopn}
\usepackage{amstext}

\DeclareMathOperator{\Lip}{Lip}
\DeclareMathOperator{\Vol}{Vol}
\DeclareMathOperator{\I}{I}
\DeclareMathOperator{\PSH}{PSH}  

\def\eps{\ensuremath{\varepsilon}}
\def\Om{\ensuremath{\Omega} }

\def\xhat{\hat{x}}
\def\Xhat{\hat{X}}
\def\muhat{\hat{\mu}}
\def\del{\partial}
\def\delbar{\bar{\partial}}

\def\R{\ensuremath{{\bf R}} }
\def\C{\ensuremath{{\bf C}} }
\def\P{\ensuremath{{\bf CP}} }

\def\M{\ensuremath{{\mathcal M}}}

\def\ddc{{dd^c}}

\theoremstyle{plain}
\newtheorem{thm}{Theorem}
\newtheorem{lem}{Lemma}
\newtheorem{conj}{Conjecture}

\theoremstyle{definition}

\theoremstyle{remark}

\begin{document}

\title[Dimension of pluriharmonic measure]{Dimension of pluriharmonic measure and polynomial
                    endomorphisms of $\C^n$.}
\author{I.~Binder}
\thanks{The first author is supported by 
N.S.F. Grant DMS-9970283.}
\address{I.~Binder\\
Department of Mathematics\\
University of Illinois, Urbana-Champaign\\
1409~W.~Green Street (MC-382) \\
Urbana, Illinois 61801, USA}
\email{ilia@math.uiuc.edu}

\author{L.~DeMarco}
\address{L.~DeMarco\\Harvard University\\
Department of Mathematics\\
One~Oxford~St~\#~325\\
 Cambridge, MA 02138, USA}
\email{demarco@math.harvard.edu}
\date{June 8, 2002}

\begin{abstract}
Let $F$ be a polynomial endomorphism of $\C^n$ which extends
holomorphically to $\P^n$.  We prove that the dimension of $\mu_F$,
the pluriharmonic measure on the boundary of the filled Julia
set of $F$, is bounded above by $2n-1$.  
\end{abstract}

\maketitle


\section{Introduction}

The {\em dimension} of a probability measure on a metric 
space is defined as the minimal Hausdorff dimension of a 
set of full measure.  In this paper, we show that
the dimension of pluriharmonic measure in $\C^n$ is bounded
above by $2n-1$ when it arises as the measure of 
maximal entropy for a regular polynomial endomorphism.

For a compact set $K$ in $\C^n$, {\em pluriharmonic measure} 
is defined as 
 $$\mu_K := \ddc G_K \wedge \cdots \wedge \ddc G_K,$$
where $G_K$ is the pluricomplex Green's function of $K$
with pole at infinity,
$d = \del + \delbar$ and $d^c = \frac{i}{2\pi}(\delbar
-\del)$.  The support of $\mu_K$ is contained in 
the Shilov boundary of $K$.  See \S\ref{background}.

Let $F:\C^n\to\C^n$ be a regular polynomial endomorphism;
i.e., one which extends holomorphically to $\P^n$. 
The filled Julia set of $F$ is the compact set
  $$K_F = \{z\in\C^n: F^m(z)\not\to\infty \mbox{ as } 
                              m\to\infty \}.$$
Pluriharmonic measure $\mu_F$ on $K_F$ is $F$-invariant and of 
maximal entropy \cite{FS}.  It is not difficult to construct
examples where the Hausdorff dimension of the support of $\mu_F$
is any value up to and including $2n$.  In answer to a question 
posed in \cite{FSsurvey}, we prove:

\begin{thm}  \label{Thm1}
The dimension of pluriharmonic measure on the filled Julia set 
of a regular polynomial endomorphism of $\C^n$ is at most $2n-1$.
\end{thm}

\medskip
In \cite{Oks}, 
Oksendal conjectured that the dimension of harmonic measure
in $\C$ would never exceed 1, though the Hausdorff dimension
of its support can take values up to and including 2. 
Makarov addressed this question for simply 
connected domains, showing that the dimension of harmonic
measure is always equal to one \cite{MakThm}. 
The theorem was extended by Jones and Wolff, establishing that 
the dimension is no greater than $1$ for general planar domains 
\cite{JW}.  Moreover, Wolff proved that there is always a set of
full harmonic measure with $\sigma$-finite Hausdorff 
$1$-measure \cite{W}. 

The complex structure on the plane plays a crucial role 
in the proof of these theorems. Namely, one relies heavily 
on the subharmonicity of the function $\log|\nabla u|$ for
harmonic $u$.

It is also possible to take a dynamical approach 
to the dimension estimates. It follows from results of 
Carleson, Jones, and Makarov
(\cite{CJ} and \cite{MakSurvey}) that any planar domain can be approximated 
in some sense by domains invariant under hyperbolic
dynamical systems (the {\em fractal approximation}).
In the special case of {\em conformal Cantor sets},
Carleson obtained dimension estimates using 
the dynamics \cite{car}.  Recently, it was shown 
that it suffices to consider polynomial Julia sets 
in the fractal approximation \cite{IP}.

Harmonic measure (evaluated at infinity) on the Julia set of a polynomial
is the unique measure of maximal entropy \cite{Brolin}, \cite{Lyubich}. 
The estimate on dimension in this case 
follows from a relation to the 
Lyapunov exponent and the entropy.  Indeed, for any 
polynomial map $F$ on $\C$, we have
\begin{equation} \label{mane}
  \dim \mu_F = \frac{\log (\deg F)}{L(F)},
\end{equation}
where $\mu_F$ denotes harmonic measure on the Julia set of $F$ and 
$L(F) = \int
\log |F'| \,d\mu_F$ is the Lyapunov exponent \cite{Mane},
\cite{man}.  The Lyapunov exponent of a polynomial is bounded 
below by $\log (\deg F)$ with equality if and only if 
the Julia set is connected \cite{prz}.  

On the other hand, for harmonic measure in $\R^n$, the 
methods applied to the study of dimension in $\C$ 
fail dramatically. 
The logarithm of the gradient of a harmonic function in
$\R^n$, $n>2$, is not subharmonic, and there is no dynamical 
interpretation of harmonic measure.  
Furthermore, in \cite{Wolff}, Wolff showed 
that for each $n>2$ there exists a domain in $\R^n$ with the 
dimension of harmonic measure strictly greater than $n-1$. 
A result of Bourgain, however, gives 
an upper bound on the dimension of harmonic measure in the form 
$n-\eps(n)$ \cite{Bourgain}.  

Because of the harmonicity of $|\nabla u|^{\frac{n-2}{n-1}}$ 
for a harmonic function $u$ in $\R^n$ (see \cite{Stein}), it 
is conjectured that the dimension of harmonic measure in $\R^n$ 
does not exceed $n-1+\frac{n-2}{n-1}$. 

For pluriharmonic measure in $\C^n$, both of the observations 
which led to proofs 
of the Oksendal conjecture in $\C$ are valid: 
the measure depends on the complex structure of $\C^n$ and 
is the measure of maximal entropy for polynomial dynamics.  
Theorem \ref{Thm1} should be the first step in the proof of the 
following conjecture:

\begin{conj} \label{conj:pluri}
The dimension of pluriharmonic measure of domains in $\C^n$ 
is at most $2n-1$.
\end{conj}

\medskip\noindent
The example of area measure on the unit sphere 
in $\C^n$ shows that this estimate is sharp.  

We believe also that one can give a precise formula for the 
dimension of pluriharmonic measure in the dynamical case,
just as in dimension one (see formula (\ref{mane}) above).     
For a diffeomorphism of a real compact manifold, such a formula 
was obtained by Ledrappier and Young in \cite{LSY}:
\begin{equation} \label{eq:LY}
     \dim\mu=\sum \frac{h_i(\mu)}{\lambda_i(\mu)},
\end{equation}
where $\mu$ is an ergodic invariant measure without zero Lyapunov 
exponents, $\lambda_i(\mu)$ are its positive Lyapunov exponents, 
and $h_i(\mu)$ are the corresponding {\em directional entropies} 
(defined in \cite{LSY}).  A similar formula was established 
by Bedford, Lyubich, and Smillie for polynomial 
diffeomorphisms of $\C^2$ \cite{BLS}.  
We make the following conjecture, which would imply Theorem
\ref{Thm1}.  

\begin{conj} \label{conj2}
For any holomorphic $F:\P^n\to\P^n$ of (algebraic) degree $d>1$,
 $$\dim \mu_F =\log d\sum_{i=1}^n\frac1{\lambda_i},$$
where $\lambda_i,\ i=1\dots n$ are the Lyapunov exponents of 
$F$ with respect to $\mu_F$ repeated with multiplicities. 
\end{conj}

\medskip\noindent
{\bf Sketch proof of Theorem \ref{Thm1}}.  We rely on 
estimates on the Lyapunov exponents of $F$ with respect
to $\mu_F$.  In particular, Briend and Duval showed that
the minimal Lyapunov exponent $\lambda_{\min}$ is bounded
below by $\frac12 \log d$ (where $d$ is the degree of $F$)
\cite{BD}.  Bedford and Jonsson proved that the sum $\Lambda$
of the Lyapunov exponents satisfies $\Lambda \geq \frac{n+1}
{2} \log d$ \cite{BJ}.  Combining these, we have 
$\Lambda \geq \lambda_0 + \frac{n-1}{2}\log d$,
where $\lambda_0 = \max\{\lambda_{\max}, \log d\}$.  

We define an invariant set $Y$ of full measure so that 
preimages of small balls centered at points in $Y$ scale
in  a way governed by the Lyapunov exponents. Namely, for each point $y\in Y$  
there exists an infinite set $M_y\subset{\mathbb Z}$ such that if $m\in M_y$
 then the $m$-th
preimage of a ball of radius $r$ centered at $F^m(y)$ should contain a ball 
of radius $\approx r e^{-m\lambda_{\max}}$ around $y$. 
In addition, the component of the preimage containing $y$ will have volume 
$\approx r^{2n} e^{-2m\Lambda}$.  
The details of the construction are very similar to the methods of \cite{BD}.

Let $A_m=\left\{y\in Y: m\in M_y\right\}$. Note that $Y=\cap_k\cup_{m\geq k}A_m$.
If we cover $Y$ by $N$ balls of radius $r$, then the ``good'' 
(as described above) $m$-th
preimages define a cover of $A_m$ by at most $N d^{mn}$ regions of controlled
shape.  Their union contains an $r e^{-m\lambda_{\max}}$-neighborhood 
of $A_m$, and has volume $\leq N r^{2n}d^m e^{2m(n-1)\lambda_0}$, by the 
estimates above.

Finally, a standard connection between the rate of decay of volume 
of a neighborhood of $Y$ and its dimension allows us
to conclude that $\dim Y\leq 2n-1$.  
\qed

\bigskip
\section{Pluriharmonic measure in $\C^n$ and dynamics}
\label{background}

In this section, we give some of the necessary background
on pluripotential theory in $\C^n$ and its relation to 
polynomial dynamics.  More details can be found in \cite{klimek}, \cite{BT} and \cite{Bsurvey}.

Let $\PSH(\C^n)$ denote the class of plurisubharmonic 
functions in $\C^n$.  
For a compact set $K$ in $\C^n$, the pluricomplex Green's
function with pole at infinity is defined as 
$$G_K(z)
  = \sup \left\{ v(z): v\in \PSH(\C^n), v\leq0
 \text{ on }K,\ v(z)\leq\log\|z\|+O(1)\text{ near }\infty\right\}.$$
See, for example, \cite{klimek} or \cite{Bsurvey}.  
If $G_K$ is continuous, the set $K$
is said to be {\em regular}.  

In contrast to the one dimensional setting, $G_K$ is not 
necessarily pluriharmonic (or even harmonic) outside of $K$.  
It is, however, {\it maximal plurisubharmonic}, meaning that
if $v$ is any plurisubharmonic function on a domain $\Om$ 
compactly contained in $\C^n-K$ with $v\leq G_K$ on $\del\Om$,
then $v\leq G_K$ on $\Om$.  Equivalently, the Monge-Ampere
mass of $G_K$, 
 $$\mu_K = (\ddc G_K)^n,$$
vanishes in $\C^n- K$.  We call the measure $\mu_K$ the
{\em pluriharmonic measure} on $K$,
and note that it is supported in the Shilov boundary of $K$.  
In fact, if $K$ is regular, then its support is equal to the 
Shilov boundary \cite{BT}.  

Pluriharmonic measure arises in the study of dynamics just 
as in the one dimensional setting. A polynomial endomorphism 
$F:\C^n\to\C^n$ is called {\it regular} if it can be extended 
holomorphically to $\P^n$.
The {\em degree} of $F$ is the degree of its polynomial
coordinate functions.  We consider only those $F$ of degree
$> 1$.  The {\em escape rate function} of $F$ is defined by 
  $$G_F(z) = \lim_{m\to\infty} \frac{1}{d^m} \log^+ \|F^m(z)\|,$$
where $d$ is the degree of $F$ and $\log^+ = \max\{\log, 0\}$.
The function $G_F$ is continuous and agrees with the pluricomplex Green's 
function for the filled Julia set $K_F = \{z\in\C^n: 
F^m(z)\not\to\infty\}$.  Fornaess and Sibony showed that
the pluriharmonic measure $\mu_F$ on $K_F$ is ergodic for $F$ and 
a measure of maximal entropy \cite{FS}.

By the Oseledec Ergodic Theorem, $F$ has $n$ Lyapunov
exponents $\lambda_{\min} \leq \cdots \leq \lambda_{\max}$ 
almost everywhere with respect to $\mu_F$ \cite{Os}.  We will
only need the existence of the minimal, maximal, and the sum 
$\Lambda$ of the Lyapunov exponents, which we can define 
as follows:
\begin{eqnarray*}
 \lambda_{\min} 
    &=& - \lim_{m\to\infty} \frac{1}{m} 
            \int \log \|(DF^m)^{-1}\| \,d\mu_F,  \\
 \lambda_{\max}
    &=&  \lim_{m\to\infty} \frac{1}{m}
            \int \log \|DF^m\| \,d\mu_F, \mbox { and } \\
 \Lambda
    &=&  \int \log |\det DF| \,d\mu_F.
\end{eqnarray*}
Briend and Duval proved that the Lyapunov exponents are
all  positive \cite{BD}; they show
\begin{equation} \label{lmin}
  \lambda_{\min} \geq \frac12 \log d,
\end{equation}
where $d$ is the degree of $F$.  
Bedford and Jonsson studied the sum of the Lyapunov exponents, 
and demonstrate that (\cite{BJ})
\begin{equation} \label{lsum}
  \Lambda \geq \frac{n+1}{2} \log d.
\end{equation}

For the proof of Theorem \ref{Thm1}, it is convenient to work in
the {\em natural extension} $(\Xhat, F)$ where $F$ is invertible
(see \cite{CFS} and \cite{BD}).  Let $P(F) = 
\bigcup_{m\geq 0} F^m(C(F))$ be the postcritical set of $F$ 
and set $X = \C^n - \bigcup_{m\geq 0} F^{-m}(P(F))$.  The 
space $(\Xhat, F)$ is the set of all bi-infinite sequences
 $$\{\xhat = (\cdots x_{-1}x_0x_1\cdots) \in 
                \prod_{-\infty}^\infty X : F(x_i) = x_{i+1}\}.$$
The map $F$ acts on $(\Xhat, F)$ by the left shift.  We define
projections $\pi_i:(\Xhat,F) \to X$ for all $i$ by $\pi_i(\xhat)
= x_i$.  Since $\mu_F$ does not charge the 
critical locus of $F$, we have $\mu_F(X)=1$. 
The measure $\mu_F$ lifts to a unique 
$F$-invariant probability
measure $\muhat$ on $(\Xhat, F)$ so that $\pi_{0*}\muhat = \mu_F$.

\bigskip
\section{Proof of the main theorem.}

In this section we give a proof of the following theorem
which clearly implies Theorem \ref{Thm1}.  

\begin{thm}
Pluriharmonic measure $\mu_F$ on the filled Julia 
set of a degree $d$ regular polynomial endomorphism 
$F:\C^n\to\C^n$ satisfies  
 $$\dim \mu_F \leq 2n - 2 + 
         \frac{\log d}{\max\{\log d, \lambda_{\max}\}},$$
where $\lambda_{\max}$ is the largest Lyapunov exponent
of $F$ with respect to $\mu_F$.  
\end{thm}

\medskip
We begin with a classical lemma.  Statements (a) and (c)
are exactly as in \cite[Lemma 2]{BD}.   
We first observe that 
there exists a constant $C(n)$ so that for any $n\times n$
matrix $A$  with $\left\|A - \I \right\|<1$, we have 
\begin{equation}\label{eq:det}
\left|\det A -1\right|
  \leq \frac{C(n)}{2} \left\|A - \I \right\|.
\end{equation}

\medskip
\begin{lem} \label{glemma}
Let $g:\Om\to\C^n$ be a function with bounded $C^2$-norm on 
a domain $\Omega\subset\C^n$ and set 
$M = C(n)\left(\|g\|_{\C^2}+1\right)$.
Let $x\in\Om$ be a noncritical point of $g$.
Given $\eps>0$, let $r(x)=\frac{1- e^{-\eps/3}}
{2M\|(D_x g)^{-1}\|^2}$, set $B_0=B(g(x), r(x))$ and let $B_1$ 
be the preimage of $B_0$ under $g$ containing $x$. Then
\begin{enumerate}
\item[(a)] $g^{-1}$ is well defined in $B_0$,
\item[(b)] $\Lip(g|B_1) \leq \|(D_x g)\| e^{\eps/3}$,
\item[(c)] $\Lip(g^{-1}|B_0) \leq \|(D_x g)^{-1}\| e^{\eps/3}$ , and
\item[(d)] $\inf_{y\in B_1} |\det(D_y g)| 
                 \geq |\det(D_x g)|e^{-\eps/3}$.
\end{enumerate} 
\end{lem}

\medskip
\begin{proof}
Let us consider a ball $B_2=B(x,\rho)$, where 
$$\rho=\frac{e^{\eps/3} - 1}{M\|(D_x g)^{-1}\|}.$$
For each $y\in B_2$, we have
\begin{equation}  \label{eq:IDD}
 \left\|\I-(D_xg)^{-1}(D_yg)\right\|
  \leq \left(\|g\|_{C^2} + 1\right) \|(D_xg)^{-1}\| \rho 
  \leq  \frac{e^{\eps/3}-1}{C(n)}, 
\end{equation}
and in particular, $\Lip \left(\I - (D_xg)^{-1}\circ g\right) < 1$.  
If $g(y_1)=g(y_2)$ for some $y_1 \not= y_2\in B_2$, then
$$ \|y_1-y_2\|=
 \left\|\left(y_1-(D_xg)^{-1}g(y_1)\right) 
        - \left(y_2-(D_xg)^{-1}g(y_2)\right)\right\|
   < \|y_1-y_2\|, $$
which is a contradiction, and therefore $g$ is injective on $B_2$.  

To establish (a), we need to know that $B_0\subset g(B_2)$. 
Map $g$ is open on $B_2$, so it is enough to check that
if $|y_1-x|=\rho$, then $|g(y_1)-g(x)|>r(x)$. But this is 
again a direct consequence of \eqref{eq:IDD}.

Now, since $B_1\subset B_2$, we have for all $y\in B_1$,
$$ \left\|D_xg-D_yg\right\| 
   \leq  \|D_xg\| \left \|\I-(D_xg)^{-1}(D_yg)\right\|
   \leq  \|D_xg\|\frac{e^{\eps/3}-1}{C(n)},  $$
and we conclude that 
$$ \|D_yg\| \leq \|D_xg\| + \frac{e^{\eps/3}-1}{C(n)}\|D_xg\|
            \leq e^{\eps/3} \|D_xg\|,$$
establishing (b).  

To prove (c),  observe that by \eqref{eq:IDD} for $y\in B_2$,
\begin{eqnarray*}
 \|(D_yg)^{-1}\| 
   &\leq&  \|(D_xg)^{-1}\| \|(\I- (\I-(D_xg)^{-1}D_yg)^{-1}\| \\
   &\leq&  \frac{\|(D_xg)^{-1}\|}{1- \|\I- (D_xg)^{-1}D_yg\|} \\
   &\leq&  \|(D_x g)^{-1}\| e^{\eps/3}.
\end{eqnarray*}
For (d), 
we compute for all $y\in B_1$ (using (\ref{eq:det}) above),
\begin{eqnarray*}
|\det D_yg - \det D_xg| 
   &=& |\det D_xg|\left|1 - \det (D_xg)^{-1}D_yg\right|  \\
   &\leq& |\det D_xg| \frac{C(n)}{2} \|\I - (D_xg)^{-1}D_yg\|  \\
   &\leq& |\det D_xg| \frac{1}{2} (e^{\eps/3}-1) \\
   &\leq& |\det D_xg| (1-e^{-\eps/3}),
\end{eqnarray*}
and therefore,
  $$\inf_{y\in B_1} |\det D_yg| \geq |\det D_xg| e^{-\eps/3}.$$
\end{proof}

\bigskip

Let $F$ be a regular polynomial endomorphism of $\C^n$ 
and $\mu_F$ the pluriharmonic measure on the boundary of the 
filled Julia set of $F$.  Denote by $\lambda_{\min}$, 
$\lambda_{\max}$, and $\Lambda$,  the minimal, maximal, 
and sum of the $n$ Lyapunov exponents of $F$ with respect
to $\mu_F$.  The space $(\Xhat, F)$ denotes the natural 
extension of $F$.  See \S\ref{background}.  

\medskip
\begin{lem}  \label{Lemma2}
Given $\eps>0$, there exist 
measurable functions $r$ and $\kappa$ on $(\Xhat,F)$ so that
for almost every point $\xhat$, we have 
$r(\xhat)>0$ and $\kappa(\xhat)<\infty$, 
and for each $m\geq 0$, a well-defined branch of 
$F^{-m}$ sending $x_0$ to $x_{-m}$ with 
\begin{enumerate}
\item[(a)]   $F^{-m}(B(x_0,s)) \supset 
        B(x_{-m},\frac{s}{\kappa(\xhat)} e^{-m(\lambda_{\max}+\eps)})$
	for all $s\leq r(\xhat)$, and 
\item[(b)]   $\Vol F^{-m}B(x_0,r(\xhat)) \leq 
        \kappa(\xhat) e^{-m(2\Lambda -\eps)}$.  
\end{enumerate}
\end{lem}

\medskip
\begin{proof}
Choose $N$ so that 
 $$0< \lambda_{\min} - \eps 
   \leq -\frac{1}{N}\int \log\|(DF^N)^{-1}\| \,d\mu_F
   \leq \lambda_{\min}, $$
and
 $$ \lambda_{\max} 
    \leq \frac{1}{N} \int \log\|DF^N\| \,d\mu_F 
    \leq \lambda_{\max} + \eps.$$
Observe that 
 $$\Lambda = \frac{1}{N} \int \log |\det DF^N| \,d\mu_F$$
for any $N\geq 0$.  

For notational simplicity, set $g = F^N$.  Observe that it is enough to prove the statement of the lemma for $g$ instead of $F$.

Fix $\xhat\in (\Xhat,g)$.  Let 
 $$r(x_{-m}) = \frac{1-e^{-\eps/3}}{2M \|(D_{x_{-m}}g)^{-1}\|^2},$$
as in Lemma \ref{glemma} where $\Omega$ is a large ball 
containing the filled Julia set of $F$.  By the ergodic theorem
applied to the function 
  $$\xhat \mapsto \log \|(D_{x_0}g)^{-1}\|,$$ 
we have 
  $$\lim_{m\to\infty} \frac{1}{m} \log \|(D_{x_{-m}}g)^{-1}\| = 0,$$
and therefore there exists a measurable function $\eta>0$ on $(\Xhat,g)$
with
\begin{equation} \label{rdecay}
  r(x_{-m}) \geq \eta(\xhat)e^{-m\eps/2}
\end{equation}
for all $m\geq 0$ and almost every $\xhat$.  

Let $B_m = B(x_0,r(x_{-1}))\cap \cdots \cap g^mB(x_{-m},r(x_{-m-1}))$.
Let $g^{-m}$ denote the inverse branch of $g^m$ taking $x_0$ 
to $x_{-m}$, well-defined on $B_m$ by Lemma \ref{glemma}(a).  
Iterating results (b), (c), 
and (d) of Lemma \ref{glemma}, we have 
 $$ \Lip(g^{-m}|B_m) \leq \|(D_{x_{-m}}g)^{-1}\| \cdots
                         \|(D_{x_{-1}}g)^{-1}\| e^{m\eps/3}, $$
 $$ \Lip(g^m|g^{-1}B_m) \leq \|D_{x_{-m}}g\| \cdots
			  \|D_{x_{-1}}g\| e^{m\eps/3}, $$
and
 $$ \inf_{y\in g^{-m}B_m} | \det D_yg^m| \geq
                        |\det D_{x_{-m}}g^m | e^{-m\eps/3}.$$

Applying the ergodic theorem to the functions 
$\xhat\mapsto \log\|(D_{x_0}g)^{-1}\|$, $\xhat\mapsto 
\log\|D_{x_0}g\|$, and $\xhat\mapsto \log|\det D_{x_0}g|$,
we see that there exists a measurable function $1 \leq C(\xhat)
< \infty$ so that 
\begin{equation} \label{g-m}
  \Lip (g^{-m}|B_m) \leq C(\xhat) e^{-m(N\lambda_{\min} - \eps/2)},
\end{equation}
\begin{equation} \label{gm}
  \Lip (g^m|g^{-m}B_m) \leq C(\xhat) e^{m(N\lambda_{\max} + \eps/2)},
\end{equation}
and
\begin{equation} \label{detg}
  \inf_{y\in g^{-m}B_m} |\det D_yg^m| \geq 
           \frac{1}{C(\xhat)} e^{m(N\Lambda - \eps/2)}  
\end{equation}
for almost every $\xhat$.  

Let $r(\xhat) = \min\{ \eta(\xhat)/C(\xhat), 1\}$.  By induction
and the estimates (\ref{rdecay}) and (\ref{g-m}), we establish that 
$B(x_0,r(\xhat))$ is contained in $B_m$ for all $m\geq 0$.  
By (\ref{gm}), we have 
 $$B(x_0,r(\xhat)) \supset g^m B(x_{-m}, \frac{r(\xhat)}{2C(\xhat)}
                            e^{-m(N\lambda_{\max}+\eps/2)}).$$
By (\ref{detg}), the volume of $g^{-m}B(x_0,r(\xhat))$ is bounded by
 $$ \Vol (g^{-m}B(x_0,r(\xhat))) 
    \leq \Vol (B(x_0,r(\xhat))) C(\xhat)^2 e^{-m(2N\Lambda - \eps)}.$$
The Lemma is proved upon setting $\kappa(\xhat) = 2C(\xhat)^2 \Vol B_1$.
\end{proof}

\bigskip
\noindent{\bf Proof of Theorem 1}.  
For fixed $\eps>0$, let $r$ and $\kappa$ be as in Lemma
\ref{Lemma2}.  Let $d$ be the degree of $F$ and   
$\lambda_{\min}$, $\lambda_{\max}$, and $\Lambda$
the minimal, maximal, and sum of the Lyapunov exponents of 
$F$.

Choose a 
set $\hat{A}\subset (\Xhat,F)$ and $r_0, \kappa_0 >0$ so that  
 $$\hat{A} \subset \{ \xhat\in (\Xhat,F): r(\xhat)\geq r_0
                    \mbox{ and } \kappa(\xhat)\leq \kappa_0\},$$
$\pi_0(\hat{A})$ has compact closure in $\C^n$ and $\muhat(\hat{A})>0$.  
Let $\hat{Y}\subset (\Xhat,F)$ be the set of all points whose 
forward orbit under $F$ lands in $\hat{A}$ infinitely often.
By ergodicity, $\muhat(\hat{Y}) = 1$.  Let $Y = \pi_0(\hat{Y}) 
= \{x_0:\xhat\in\hat{Y}\}$, so $\mu(Y) = 1$, and let 
$A_i = \pi_{-i}(\hat{A})$.  Observe that 
  $$Y = \bigcap_{l\geq 0} \bigcup_{m\geq l} A_m.$$

We will show that the Hausdorff dimension of $Y$ is bounded 
above by $2n-2 + \frac{\log d}{\lambda_0} + \frac{4\eps}{\lambda_0}$, 
where
$\lambda_0 = \max\{\lambda_{\max}, \log d\}$.  As $Y$ has 
full measure and $\eps$ is arbitrary, this will prove the 
theorem.  

For a ball $B$ in $\C^n$, let $\frac12 B$ denote a concentric
ball with half the radius.  Let $\Sigma$ denote a finite 
collection of balls $B$ of radius $r_0$ so that the balls 
$\frac12 B$ cover 
$A_0$.  For each point $y\in A_m$, select $\hat{y}\in\hat{A}$
so that $y=\pi_{-m}(\hat{y})$.  Choose element $B$ of $\Sigma$
so that $\pi_0(\hat{y})$ lies in $\frac12 B$.  Let $B_y$
be the preimage of $F^{-m}B$ containing $y$.  The collection 
of these $B_y$ for all $y\in A_m$ defines the finite 
cover $\Sigma_m$ of $A_m$.  

If $\sigma$ is the number of elements in $\Sigma$, then 
the number of elements of $\Sigma_m$ is no greater than
$\sigma d^{mn}$.  Let $\lambda_0 = \max\{\log d, \lambda_{\max}\}$.
We will establish the following two properties of the 
cover $\Sigma_m$:  
\begin{enumerate}
\item[(I)]  The union $\bigcup_{B\in\Sigma_m} B$ contains an 
            $\frac{r_0}{4\kappa_0} e^{-m(\lambda_0+\eps)}$-neighborhood
            of $A_m$, and 
\item[(II)] $\Vol \left(\bigcup_{B\in\Sigma_m} B\right) \leq
            \sigma d^m \kappa_0 e^{-m(2\lambda_0 - \eps)}$.
\end{enumerate}
Observe first that for each $y\in A_m$, the set $B_y\in\Sigma_m$
contains a ball of radius $\frac{r_0}{4\kappa_0}
e^{-m(\lambda_{\max}+\eps)}$ around $y$ by Lemma \ref{Lemma2}.  
Of course, $\lambda_{\max}\leq\lambda_0$, giving (I).  

To establish (II), we observe that as $F^mB_y\subset 
B(\pi_0(\hat{y}),r_0)$ for each $B_y\in\Sigma_m$, 
Lemma \ref{Lemma2}(b) implies that 
  $$\Vol B_y \leq \kappa_0 e^{-m(2\Lambda-\eps)}.$$
Summing over the volumes of all elements in $\Sigma_m$,
we write 
  $$\Vol \left(\bigcup_{B\in\Sigma_m} B\right) 
        \leq \sigma d^{mn} \kappa_0 e^{-m(2\Lambda -\eps)}.$$
By (\ref{lsum}), $\Lambda$ is bounded below by $\frac{n+1}{2} 
\log d$, and by (\ref{lmin}), each Lyapunov exponent is bounded
below by $\frac12 \log d$.  Combining these gives 
$\Lambda\geq \frac{n-1}{2} \log d + \lambda_0$ and we obtain
statement (II).  

We define a covering $\M_m$ of $A_m$ to be the collection of 
all mesh cubes of edge length $\frac{1}{\sqrt{2n}}
\frac{r_0}{4\kappa_0} e^{-m(\lambda_0+\eps)}$ which intersect
$A_m$.  Let $c = \frac{1}{\sqrt{2n}}\frac{r_0}{4\kappa_0}$.
By property (I), each cube is contained in an element 
of $\Sigma_m$.  The number of cubes in $\M_m$ is bounded
above by the volume $\Vol \left(\cup_{B\in\Sigma_m} B\right)$
divided by the volume of each cube.  That is,
\begin{eqnarray*}
 |\M_m| 
 &\leq& \frac{\sigma d^m \kappa_0 e^{-m(2\lambda_0-\eps)}}
           {c^{2n} e^{-2mn(\lambda_0 +\eps)}}   \\
 &=&   \frac{\sigma\kappa_0}{c^{2n}} d^m e^{2(n-1)m\lambda_0}
	                  e^{(2n+1)m\eps}. 
\end{eqnarray*}

We now show that Hausdorff measure of $Y$ in dimension 
$2n-2 + \frac{\log d}{\lambda_0} + \frac{4\eps}{\lambda_0}$ is finite,
completing the proof.
Fix $\delta>0$.  Choose $l\geq 0$ so that the mesh cubes
in $\M_m$ are of diameter $\delta_m \leq \delta$ for each $m\geq l$. 
The union of the elements of $\M_m$ for $m\geq l$ covers $Y$.
Therefore, 
\begin{eqnarray*}
 m_{2n-2 + \frac{\log d}{\lambda_0} + \frac{4\eps}{\lambda_0}}(Y)
 &\leq&  \sum_{m\geq l} |\M_m| 
    (\delta_m)^{2n-2 + \frac{\log d}{\lambda_0} + \frac{4\eps}{\lambda_0}} \\
 &\leq&   C \sum_{m\geq l} d^m e^{2(n-1)m\lambda_0} e^{(2n+1)m\eps}
                  e^{-m(\lambda_0+\eps)\left(2n-2 + 
          \frac{\log d}{\lambda_0} + \frac{4\eps}{\lambda_0}\right)}   \\
 &=&  C \sum_{m\geq l} e^{-m\eps\left(1 + \frac{\log d}{\lambda_0}
                     +   \frac{4\eps}{\lambda_0}\right) }  \\
 &\leq&  C \sum_0^\infty e^{-m\eps} < \infty
\end{eqnarray*}
\qed

\bigskip

\providecommand{\bysame}{\leavevmode\hbox to3em{\hrulefill}\thinspace}
\providecommand{\MR}{\relax\ifhmode\unskip\space\fi MR }
\providecommand{\MRhref}[2]{%
  \href{http://www.ams.org/mathscinet-getitem?mr=#1}{#2}
}
\providecommand{\href}[2]{#2}


\begin{thebibliography}{{O}ks81}

\bibitem[BD99]{BD}
Jean-Yves Briend and Julien Duval, \emph{Exposants de {L}iapounoff et
  distribution des points p\'eriodiques d'un endomorphisme de $\bold {C}{ \rm
  {P}}\sp k$}, Acta Math. \textbf{182} (1999), no.~2, 143--157.
  \MR{2000f:32023}

\bibitem[Bed93]{Bsurvey}
Eric Bedford, \emph{Survey of pluri-potential theory}, Several complex
  variables (Stockholm, 1987/1988), Princeton Univ. Press, Princeton, NJ, 1993,
  pp.~48--97. \MR{94b:32014}

\bibitem[BJ00]{BJ}
Eric Bedford and Mattias Jonsson, \emph{Dynamics of regular polynomial
  endomorphisms of $\bf {C}\sp k$}, Amer. J. Math. \textbf{122} (2000), no.~1,
  153--212. \MR{2001c:32012}

\bibitem[BJ02]{IP}
I.~Binder and P.~Jones, In preparation, 2002.

\bibitem[BLS93]{BLS}
Eric Bedford, Mikhail Lyubich, and John Smillie, \emph{Polynomial
  diffeomorphisms of ${\bf {c}}\sp 2$. {I}{V}. {T}he measure of maximal entropy
  and laminar currents}, Invent. Math. \textbf{112} (1993), no.~1, 77--125.
  \MR{94g:32035}

\bibitem[Bou87]{Bourgain}
J.~Bourgain, \emph{On the {H}ausdorff dimension of harmonic measure in higher
  dimension}, Invent. Math. \textbf{87} (1987), no.~3, 477--483. \MR{88b:31004}

\bibitem[Bro65]{Brolin}
Hans Brolin, \emph{Invariant sets under iteration of rational functions}, Ark.
  Mat. \textbf{6} (1965), 103--144 (1965). \MR{33 \#2805}

\bibitem[BT87]{BT}
Eric Bedford and B.~A. Taylor, \emph{Fine topology, \v {S}ilov boundary, and
  $(dd\sp c)\sp n$}, J. Funct. Anal. \textbf{72} (1987), no.~2, 225--251.
  \MR{88g:32033}

\bibitem[Car85]{car}
L.~Carleson, \emph{On the support of harmonic measure for the sets of {C}antor
  type}, Ann. Acad. Sci. Fenn \textbf{10} (1985), 113--123.

\bibitem[CFS82]{CFS}
I.~P. Cornfeld, S.~V. Fomin, and Ya.~G. Sina{\u\i}, \emph{Ergodic theory},
  Springer-Verlag, New York, 1982, Translated from the Russian by A. B.
  Sosinski\u\i. \MR{87f:28019}

\bibitem[CJ92]{CJ}
Lennart Carleson and Peter~W. Jones, \emph{On coefficient problems for
  univalent functions and conformal dimension}, Duke Math. J. \textbf{66}
  (1992), no.~2, 169--206. \MR{93c:30022}

\bibitem[FS95]{FS}
John~Erik Fornaess and Nessim Sibony, \emph{Complex dynamics in higher
  dimension. {I}{I}}, Modern methods in complex analysis (Princeton, NJ, 1992),
  Princeton Univ. Press, Princeton, NJ, 1995, pp.~135--182. \MR{97g:32033}

\bibitem[FS01]{FSsurvey}
\bysame, \emph{Complex dynamics in higher dimension}, Preprint, 2001.

\bibitem[JW88]{JW}
Peter~W. Jones and Thomas~H. Wolff, \emph{Hausdorff dimension of harmonic
  measures in the plane}, Acta Math. \textbf{161} (1988), no.~1-2, 131--144.

\bibitem[Kli91]{klimek}
Maciej Klimek, \emph{Pluripotential theory}, The Clarendon Press Oxford
  University Press, New York, 1991, Oxford Science Publications. \MR{93h:32021}

\bibitem[Lju83]{Lyubich}
M.~Ju. Ljubich, \emph{Entropy properties of rational endomorphisms of the
  {R}iemann sphere}, Ergodic Theory Dynam. Systems \textbf{3} (1983), no.~3,
  351--385. \MR{85k:58049}

\bibitem[LY85]{LSY}
F.~Ledrappier and L.-S. Young, \emph{The metric entropy of diffeomorphisms.
  {I}{I}. {R}elations between entropy, exponents and dimension}, Ann. of Math.
  (2) \textbf{122} (1985), no.~3, 540--574. \MR{87i:58101b}

\bibitem[Mak85]{MakThm}
N.~G. Makarov, \emph{On the distortion of boundary sets under conformal
  mappings}, Proc. London Math. Soc. (3) \textbf{51} (1985), no.~2, 369--384.

\bibitem[Mak89]{MakSurvey}
\bysame, \emph{Probability methods in the theory of conformal mappings.},
  Algebra i Analiz \textbf{1} (1989), no.~1, 3--59.

\bibitem[Man84]{man}
Anthony Manning, \emph{The dimension of the maximal measure for a polynomial
  map}, Ann. of Math. (2) \textbf{119} (1984), no.~2, 425--430. \MR{85i:58068}

\bibitem[Mn88]{Mane}
Ricardo Ma\~{n}{\'e}, \emph{The {H}ausdorff dimension of invariant
  probabilities of rational maps}, Dynamical systems, Valparaiso 1986,
  Springer, Berlin, 1988, pp.~86--117. \MR{90j:58073}

\bibitem[{O}ks81]{Oks}
Bernt {O}ksendal, \emph{Brownian motion and sets of harmonic measure zero},
  Pacific J. Math. \textbf{95} (1981), no.~1, 179--192. \MR{83c:60106}

\bibitem[Ose68]{Os}
V.~I. Oseledec, \emph{A multiplicative ergodic theorem. {C}haracteristic
  {L}japunov, exponents of dynamical systems}, Trudy Moskov. Mat. Ob\v s\v c.
  \textbf{19} (1968), 179--210. \MR{39 \#1629}

\bibitem[Prz85]{prz}
Feliks Przytycki, \emph{Hausdorff dimension of harmonic measure on the boundary
  of an attractive basin for a holomorphic map}, Invent. Math. \textbf{80}
  (1985), no.~1, 161--179. \MR{86g:30035}

\bibitem[Ste70]{Stein}
Elias~M. Stein, \emph{Singular integrals and differentiability properties of
  functions}, Princeton University Press, Princeton, N.J., 1970. \MR{44 \#7280}

\bibitem[Wol93]{W}
Thomas~H. Wolff, \emph{Plane harmonic measures live on sets of $\sigma$-finite
  length}, Ark. Mat. \textbf{31} (1993), no.~1, 137--172.

\bibitem[Wol95]{Wolff}
\bysame, \emph{Counterexamples with harmonic gradients in ${{\bf {{R}}}}\sp
  3$}, Essays on Fourier analysis in honor of Elias M. Stein (Princeton, NJ,
  1991), Princeton Univ. Press, Princeton, NJ, 1995, pp.~321--384.
  \MR{95m:31010}

\end{thebibliography}
\end{document}